\documentclass[a4paper,twoside,final]{article}
\usepackage{amssymb}
\usepackage{amsmath}
\usepackage{amsthm}
\usepackage{amscd}
\usepackage[utf8]{inputenc}
\usepackage[T1]{fontenc}
\usepackage{lpretex}
\usepackage[all]{xy}
\usepackage{title}

\newcommand\Kscheme[2][]{K_0^{#1}(\mathrm{Spc}_{#2})}
\newcommand\Kschemecompl[2][]{\widehat{K}_0^{#1}(\mathrm{Spc}_{#2})}
\newcommand\Bclass{\mathrm{B}}
\newcommand\Eclass{\mathrm{E}}
\newcommand\Lclass{\mathbb{L}}
\newcommand\coeffs{\widehat\Z'}
\begin{document}
\title[Approximating classifying spaces]
{Approximating classifying spaces\\ by\\ smooth projective varieties}

\author{Torsten Ekedahl} 

\address{Matematiska institutionen\\
Stockholms universitet\\ 
SE-106 91 Stockholm\\ 
Sweden}

\email{teke@math.su.se} 
%% 14F45 Topological properties
%% 14F35 Homotopy theory; fundamental groups
%% 14F25 Classical real and complex cohomology
%% 14L24 Geometric invariant theory [See also 13A50]
%% 14L30 Group actions on varieties or schemes (quotients) 
%% 55R35 Classifying spaces of groups and ${H}$-spaces
%% 55R40 Homology of classifying spaces, characteristic classes

\begin{abstract}
We prove that for every reductive algebraic group $H$ with centre of positive
dimension and every integer $K$ there is a smooth and projective variety $X$ and
an algebraic $H$-torsor $P \to X$ such that the classifying map $X \to \Bclass
H$ induces an isomorphism in cohomology in degrees $\le K$. This is then applied
to show that if $G$ is a connected non-special group there is a $G$-torsor $P
\to X$ for which we do not have $[P]=[G][X]$ in the (completion of the)
Grothendieck ring of varieties.
\end{abstract}
\maketitle

For a fixed integer $k$ and a variable integer $n$, the Grassmannian
$\text{Gr}(k,n+k)$ of (complex) $k$-spaces in (complex) $n+k$-dimensional space
becomes a closer and closer (the larger we make $n$) homotopy approximation to
the classifying space of $\GL_k(\C)$ in the following sense: The frame bundle of
the tautological $k$-bundle is a $\GL_k(\C)$-torsor over $\text{Gr}(k,n+k)$ and
thus gives a map $\text{Gr}(k,n+k) \to \Bclass\GL_k(\C)$ which induces an
isomorphism on homotopy groups up to a degree which tends to infinity with
$n$. The point here is that $\text{Gr}(k,n+k)$ is a complex projective manifold
and that the $\GL_k(\C)$-torsor is also algebraic. Hence, we can homotopically
approximate the classifying space of $\GL_k(\C)$ by projective manifolds using
algebraic torsors. This is somewhat special to $\GL_k(\C)$; it is for instance
not possible to do it for $\SL_k(\C)$ as $H^2(\Bclass\SL_k(\C),\Z)=0$ and the
second cohomology group of a projective manifold is never zero. When $G$ is a
finite group this was rectified by Atiyah and Hirzebruch (cf.,
\cite[Prop.~6.6]{atiyah62::analy}) who showed that the Godeaux-Serre method
(cf., \cite{serre58::sur}) together with the Lefschetz hyperplane theorem could
be used to prove that we can get homotopy approximations to
$\Bclass(G\times\mul)$ which are given by algebraic torsors over a smooth and
projective base. The main result of the present paper is an extension of this
result to the case when $G$ is a reductive group (and where we also allow
ourselves to replace $G\times\mul$ by any central extension of $G$ by $\mul$).

The technical problems we have to fight with is the fact that when $\dim G>0$
there will be unstable points for the action of $G$ on a projective space. The
unstable points will contribute to the cohomology of a (suitably iterated)
hypersurface section defined by a $G$-invariant form which will have to be
``subtracted off'' if we are to get hold of the cohomology of the torsor (and
consequently its smooth projective base). This requires us to get information on
the cohomology of the unstable locus and we shall use the inductive analysis of
the unstable points introduced by Kempf et al. We are much helped by the fact
that we can allow ourselves to consider $\P(V^n)$ where $V$ is a fixed linear
representation and $n$ is arbitrarily large; making $n$ large enough means that
big pieces of the cohomology of the unstable locus will not interfere with our
goal. On the other hand, in the argument we shall use (repeatedly) Poincaré
duality and consequently our results will only give a homological and not a
homotopy approximation.

The fact that our constructions are very non-canonical makes it unlikely that
they could be used to study the the classifying spaces of reductive (or
equivalently compact) groups. Instead it can be used to show the existence of
projective manifolds exhibiting various topological behaviours (as was already
demonstrated by Atiyah and Hirzebruch where a particular behaviour of cohomology
under the action of the Steenrod algebra was required). Indeed, we shall apply
it to show that if $G$ is a linear connected non-special (in the sense of
Grothendieck and Serre) group, then there are algebraic $G$-torsors with a
smooth projective base such that the associated $G/B$-fibration has cohomology
which is not (additively) isomorphic to the cohomology of the product of the
base and $G/B$. (On the contrary such an isomorphism always exists for a special
group.) This will then be used to obtain the result on the Grothendieck ring of
varieties mentioned in the abstract.

%;;%Section The approximation
\begin{section}{The approximation}

So far we have formulated our result in topological terms and that forces us to
work with algebraic groups over the complex numbers. In Subsection
\ref{ss:Algebraic version} we shall discuss the necessary modifications to
obtain results valid over any field.

For the proof of our first proposition we need to recall (cf., \cite[Part
II]{kirwan84::cohom} whose notation we shall follow) the inductive analysis of
unstable points of a representation $W$ of a reductive group $G$. Hence, if we
fix a maximal torus $T$ of $G$ we have a stratification $\{S_\beta\}$ of $\P(W)$
parametrised by certain elements of $M\Tensor\Q$, where $M$ is the cocharacter
group of $T$. To such a $\beta\ne 0$ the following data can be associated:
\begin{itemize}
\item A positive definite quadratic form $q$ on $M\Tensor\Q$ that is invariant
under the Weyl group has been chosen once and for all.

\item Putting $N_\beta:=\set{\alpha\in N}{\alpha\cdot\beta=q(\beta)}$ and
$N'_\beta:=\set{\alpha\in N}{\alpha\cdot\beta>q(\beta)}$, where $N$ is the
character group of $T$, we let $W_\beta \subseteq W$ be the subspace spanned
by the $T$-eigenvectors of weights belonging to $N_\beta$.

\item We let $T_\beta$ be the image of $r\beta$, where $r$ is a strictly
positive integer such that $r\beta\in M$. The centraliser $\text{Stab}_\beta$ of
$T_\beta$ is a reductive group stabilising $W_\beta$. There is a reductive
subgroup $G_\beta$ of $\text{Stab}_\beta$ such that $\text{Stab}_\beta=T_\beta
G_\beta$ and $T_\beta$ is not contained in $G_\beta$. In particular
$\dim\text{Stab}_\beta=\dim G_\beta+1$.

\item Setting $Z_\beta:=\P(W_\beta) \subseteq \P(W)$ we have that
$\text{Stab}_\beta$ acts (linearly) on it. We let $Z^{ss}_\beta\subseteq
Z_\beta$ be the semi-stable locus for the action of $G_\beta$. 

\item Letting $W'_\beta\subseteq W$ be the subspace spanned by the
$T$-eigenvectors of weights belonging to $N'_\beta$, $W'_\beta+W_\beta\subseteq
W$ is a direct sum. Putting $Y_\beta :=\P(W_\beta\Dsum W'_\beta)\setminus
\P(W'_\beta)$ we have a linear projection morphism $Y_\beta \to Z_\beta$ which is
a vector bundle of constant rank $r_\beta(V):=\dim W'_\beta$. (Under the
identification of $Y_\beta \to Z_\beta$ with a vector bundle $Z_\beta$ becomes
the zero section.) We let $Y^{ss}_\beta$ be the inverse image of $Z^{ss}_\beta$
under the map $Y_\beta \to Z_\beta$.

\item There is a parabolic subgroup $P_\beta$ stabilising $Y_\beta$ and
$Y^{ss}_\beta$ as well as containing $\text{Stab}_\beta$. More precisely
$\text{Stab}_\beta$ is a Levi factor of $P_\beta$, i.e., it maps bijectively to
$P_\beta/U_\beta$, where $U_\beta$ is the unipotent radical of $P_\beta$. In
particular, as $U_\beta$ and $G/P_\beta$ have the same dimension we have that
$\dim G=\dim \text{Stab}_\beta+2\dim(G/P_\beta)$.

\item We have an isomorphism $S_\beta\iso G\times_{P_\beta}Y^{ss}_\beta$ and
$S_0$ equals $\P(W)^{ss}$, the $G$-semistable locus of $\P(W)$. In particular if
there is no weight of $W$ that lies in $N_\beta$ we have that
$S_\beta=\emptyset$.
\end{itemize}
We now denote by $S'_\beta$, $Y'^{ss}_\beta$ and $Z'^{ss}_\beta$ the inverse
images under the quotient map $W\setminus\{0\} \to \P(W)$ of $S_\beta$,
$Y^{ss}_\beta$ and $Z^{ss}_\beta$ respectively. The linear projection map
$Y_\beta \to Z_\beta$ induces a vector bundle map $Y'^{ss}_\beta \to
Z'^{ss}_\beta$ of rank $r_\beta(W)$. Consequently we have that
%;;%%Equation YZ shift
\begin{equation}\label{YZ shift}
H^i_c(Y'^{ss}_\beta,\Z)=H^{i-2r_\beta(W)}_c(Z'^{ss},\Z)
\end{equation}
and in particular $H^i_c(Y'^{ss}_\beta,\Z)=0$ whenever $i < 2r_\beta(W)$. We also have
that $S'_\beta=G\times_{P_\beta}Y'_\beta$.  Its pullback along $G \to
G/P_\beta$ is a product and as $P_\beta$ is connected, the proper base change
theorem implies that $R^i\pi_!\Z$ is constant with constant value
$H^i_c(Y'^{ss}_\beta,\Z)$. This (together with the fact that $G/P_\beta$ is
compact) implies that the Leray spectral sequence has the form
%;;%%Equation parabolic s.s.
\begin{equation}\label{parabolic s.s.}
H^j(G/P_\beta,H^i_c(Y'^{ss}_\beta,\Z)) \implies H^{i+j}_c(S'_\beta,\Z).
\end{equation}

Turning now to a slightly different topic let $V$ be a $k$-dimensional complex
vector space for $k>0$. We define the open subset $U_n \subseteq V^n$ by the
condition that $(v_1,\dots,v_n) \in U_n$ if there is a subset $S\subseteq
\{1,\dots,n\}$ with $|S|=k$ and such that $(v_s)_{s\in S}$ is a basis for
$V$. Clearly, $U_n$ is stable under multiplication by scalars and so corresponds
to an open subset $W_n$ of $\P(V^n)$.
%;;%%Lemma Stability (iv) Stable topology
\begin{lemma}\label{Stability}
\part[0] A tuple $(v_1,\dots,v_n)$ lies in $U_n$ precisely when the $v_i$ span
$V$.

\part[i] $\GL(V)$ acts freely on $U_n$ and hence so does any algebraic subgroup of
$\GL(V)$.

\part[ii] Any point of $W_n$ is semi-stable (and therefore also stable by
\DHRefpart{i}) for the action by $\SL(V)$ and thus also for any algebraic subgroup
of $\SL(V)$.

\part[iii] The codimension of $V^n\setminus U_n$ in $V^n$ is $\ge n+1-k$.

\part[iv]\label{Stable topology} Let $G \subseteq \SL(V)$ be a (not necessarily
connected) reductive algebraic subgroup. Let $V_{ss}^n\subseteq
V^n\setminus\{0\}$ be the inverse image of the locus of $G$-semistable points of
$\P(V^n)$. Then for every $K$ there is an integer $N$ such that
$H^i_c(V_{ss}^n,\Z)=0$ for all $\dim G +1< i\le K$ and $n\ge N$.
\begin{proof}
To begin with \DHrefpart{0} is just linear algebra.

For \DHrefpart{i} a point $(v_1,\dots,v_n)$ of $U_n$ contains a basis and hence
its stabiliser is reduced to the identity. Hence it remains to show that the
morphism $\GL(V)\times U_n \to U_n\times U_n$ given by $(g,x) \mapsto (x,gx)$ is
proper. For this we use the valuation criterion so we assume that we have a
discrete valuation ring $R$ with fraction field $K$ and $K$-points $g$ and $x$
such that $x$ and $gx$ are defined as $R$-points of $U_n$. As $x$ and $gx$ thus
both span $V\Tensor R$ we get that $g$ and its inverse are defined over $R$.

As for \DHrefpart{ii}, for every $S \subseteq \{1,\dots,n\}$ we have the
$k\times S$-determinant which is $\SL(V)$-invariant and for every point of
$U_n$ there is an $S$ such that it separates the point from the origin.

A point $(v_1,\dots,v_n)$ in the complement of $U_n$ has, by \DHrefpart{0}, the
property that its linear span is a proper subspace of $V$. Hence, the point lies
in the union of the $U^n$, where $U$ runs over the hyperplanes of $V$. That
union has dimension at most $k-1+n(k-1)$ and hence codimension at least
$nk-(n+1)(k-1)=n+1-k$ which proves \DHrefpart{iii}.

Finally, to prove \DHrefpart{iv} we start by noticing that if $G^o$ is the group
of connected components of $G$, then a point of $\P(V^n)$ is semi-stable for $G$
precisely when it is so for $G^o$. Hence we may, to simplify, assume that $G$ is
connected.

We now apply the analysis recalled above of unstable points of the
$G$-representation $V^n$. The crucial part is that by construction
$(V^n)_\beta=(V_\beta)^n$ (as $N_\beta$ does not depend on the
representation). Furthermore, again by construction, the set of $\beta$ for
which the $S_\beta$ are non-empty is contained in a finite set which only
depends on the weights of $T$ appearing in the $V^n$ and is thus independent of
$n$. We now prove the result by induction on $\dim G$. The base case is when
$G=\{e\}$ in which case $V^n_{ss}=V^n\setminus\{0\}$ and thus
$H^i_c(V^n_{ss},\Z)=0$ unless $i=1$ or $i=2n\dim V$ which gives the result in
that case. In the general case, when $\dim G>0$, we have that $\dim G_\beta<\dim
G$ for all $\beta\ne 0$ so we may assume that the theorem is true for
$(G_\beta,V_\beta)$. The complement $S'$ of $V^n_{ss}$ in $V^n\setminus\{0\}$ is
the union of the $S_\beta$ where $\beta\ne 0$ runs over a finite set,
independent of $n$. Using the long exact sequences of cohomology with compact
supports we can analyse the vanishing of $H^i_c(S',\Z)$ in terms of the
vanishing of the $H^i(S'_\beta)$. We have that $r_\beta(V^n)=nr_\beta(V)$ so if
$r_\beta(V) \ne 0$ we get from (\ref{YZ shift}) (and (\ref{parabolic s.s.})) 
that by making $n$ large enough we may assume that $H^i_c(S'_\beta,\Z)=0$ when
$i\le K$. For those $\beta$ with $r_\beta(V)=0$ we have that
$Y'^{ss}_\beta=Z'^{ss}_\beta$. Hence (\ref{parabolic s.s.}) and the induction
assumption implies that by making $n$ large enough we may assume that
$H^i_c(S'_\beta,\Z)=0$ whenever $2\dim(G/P_\beta)+\dim G_\alpha+1<i\le K$ and
as, which was noted above, $2\dim(G/P_\beta)+\dim G_\alpha+1=\dim G$ we get that
$H^i_c(S'_\beta,\Z)=0$ when $\dim G<i\le K$. From this (and the long of exact
sequences of cohomology) we get that $H^i_c(S',\Z)=0$ for $\dim G<i\le K$. We
now have a long exact sequence
\begin{displaymath}
\cdots \to  H^{i-1}_c(V^n\setminus\{0\},\Z)\to  H^{i-1}_c(S',\Z)\to
H^{i}_c(V^n_{ss},\Z)\to H^{i}_c(V^n\setminus\{0\},\Z)\to\cdots.
\end{displaymath}
From it, and the fact that $H^i_c(V^n\setminus\{0\},\Z)=0$ for $1<i<n\dim V$ it
follows that $H^i(V^n_{ss},\Z)=0$ if $\dim G+1<i\le K$ (and $n$
is large enough).
\end{proof} 
\end{lemma}
The following result is no doubt well-known but lacking a reference we provide a
proof.
%;;%%Lemma Central representation
\begin{lemma}\label{Central representation}
Let $G$ be reductive group provided with a centrally embedded
$\mul\hookrightarrow G$. Then $G$ has a faithful linear representation $V$ such
that the central $\mul$ acts by the identity map into the scalar linear maps.
\begin{proof}
We first prove the following: Let $\sL$ be a line bundle on an affine scheme $X$
(over some affine base) and let $R$ be the ring generated (over the base) by the
sections of $\sL$ considered as functions on the associated $\mul$-torsor
$T$. Then the induced map $T \to \Spec R$ is an open embedding. Indeed, assume
first that $\sL=R t$. In that case $T$ maps into the subscheme $\Spec R[t^{-1}]$
but $R[t^{-1}]$ is equal to the affine coordinate ring of $T$ so that $T \to
\Spec R[t^{-1}]$ is an isomorphism. As $\sL$ is locally trivial this implies
that there is an open covering $\{U_i\}$ of $T$ such that $U_i \hookrightarrow X
\to \Spec R$ is an open embedding. From this if follows that it remains prove
that $X \to \Spec R$ is injective on points. That the elements of $\sL$
separates points of $\sL$ is clear however.

In the case of $\mul \hookrightarrow G$ we get an $\sL$ by considering the
functions on $G$ that are of weight $1$ with respect to the $\mul$-action. It is
a line bundle over $G/\mul$. As $\mul$ is central, $\sL$ is $G$-invariant and is
hence the union of finite dimensional subrepresentations. Applying what was just
shown we get that a large enough such subrepresentation is faithful. This
concludes our proof as by construction $\mul$ acts on it by the identity map
into scalars.
\end{proof}
\end{lemma}
Before passing to our main theorem we want to emphasise that we do not assume
that a reductive group is connected.
%;;%%Theorem Approximation
\begin{theorem}\label{Approximation}
Let $H$ be a reductive algebraic group with positive-dimensional centre over the
complex numbers. Then for every integer $K$ there is a smooth and projective
variety $X$ and an algebraic $H$-torsor $P \to X$ such that its classifying map
$X \to \Bclass H$ induces an isomorphism $H^i(\Bclass H,\Z) \to H^i(X,\Z)$ for
all $i \le K$.

If $H$ has the form $G\times\mul$, then it is possible to choose $P$ such that
the line bundle associated to the the $\mul$-torsor $P/G \to X$ is ample.
\begin{proof}
As the centre of $H$ is positive dimensional and as it is reductive, the
connected component of the centre is a torus so that we may find a copy of
$\mul$ contained in the centre. By Lemma \ref{Central representation} we may
find a faithful linear representation $V$ of $G$ which restricts to
multiplication by scalars on the $\mul$. We now put $G:=H\cap\SL(V)$. After
possibly replacing $V$ by $V^n$ we may, by Lemma \ref{Stability} and the
condition that $H$ acts faithfully on $V$, assume that there is an open subset
$U\subseteq V$, whose complement has arbitrarily high codimension, stable under
$H$ such that $H$ stabilises $U$ and acts freely on it. Furthermore, again by
the lemma, we may also assume that $U':=U/\mul \subseteq \P(V)$ consists of
$G$-stable points so that $U$ is a $H$-torsor over an open subset $U''$ of the
GIT quotient $\P(V)//G$. Finally, we may choose a suitable $k$ such that the
base point locus of $(S^kV)^G$ in $\P(V)$ is equal to the set $T$ of unstable
points and that $(S^kV)^G$ induces an embedding of the $\P(V)//G$ of $\P(V)$ by
$G$ into $\P((S^kV)^G)$.

We now construct a sequence $\P(V)=X_0\supseteq X_1 \supseteq \cdots \supseteq X_a$ of
closed subvarieties such that
\begin{itemize}
\item the codimension of $X_j$ in $\P(V)$ at any of its points outside of $T$ is
equal to $i$,

\item $X_{j+1}$ is the zero-set in $X_j$ of some $f_j\in (S^kV)^G$,

\item $X_j$ is smooth outside $T$ and

\item $X_a\setminus T$ lies inside $U'$.
\end{itemize}
That such a sequence exists is clear as the general hyperplane section of a
linear system is smooth outside the base points of the system. Furthermore, as
the codimension of $V\setminus U$ in $V$ can be made arbitrarily large, once
again by the lemma, so can $L$, the dimension of $X_a\setminus T$. We may also
by (\ref{Stable topology}) assume that $H^i_c(V\setminus T',\Z)=0$ when $\dim
H=\dim G+1<i<L$.

Now, by the fact that $X_j\setminus X_{j+1}$ is smooth and affine we get that
$H^i_c(X_j,\Z) \to H^i_c(X_{j+1},\Z)$ is an isomorphism for $i < \dim
(X_{j+1}\setminus T)$. Using the long exact sequence of cohomology (and the
$5$-lemma) this implies that $H^i_c(X_j\setminus T,\Z) \to
H^i_c(X_{j+1}\setminus T,\Z)$ is an isomorphism for $i < \dim (X_{j+1}\setminus
T)$. Hence we get that $H^i_c(\P(V)\setminus T,\Z) \to H^i_c(X_a\setminus T,\Z)$
is an isomorphism for $i<L$ which in turn implies that we have an isomorphism
$H^i_c(V\setminus T',\Z) \to H^i_c(U'_a,\Z)$, where $U'_a$ is the inverse image
of $X_a\setminus T$ in $V\setminus\{0\}$ and $T'$ the inverse image of $T$. Thus
we get that $H^i_c(U'_a)=0$ for $\dim H<i<L$. Now, by construction $H$
acts freely on $U'_a$ and the quotient $Y_a$ is the zero set in $\P(V)//G$ of
the $f_i$. Hence $Y_a$ is proper as it is closed in $\P(V)//G$ and smooth as
$U'_a$ is. Furthermore, the $Y_j$ are general iterated hyperplane sections and
hence normal (by \cite{seidenberg50} and the fact that $\P(V)//G$ is) and
therefore they are connected by the Lefschetz theorem (and of course the fact
that $\P(V)//G$ is connected). This may be applied to
the case when $H$ is replaced by its connected component $H^o$ and thus also
$U'_a$ is connected (as $U'_a \to U'_a/H^o$ is a fibration with connected
fibres). Furthermore, if $K \subseteq H$ is a maximal compact subgroup, then the
$G$-torsor $U'_a \to Y_a$ has a reduction to a $K$-torsor $Q_a\subseteq U'_a$
and this inclusion is a homotopy equivalence. Furthermore, as $Y_a$ is a compact
(complex) manifold, $Q_a$ is also a compact manifold. Now, letting $K^o:=H^o\cap
K$ the fibration $Q_a \to Q_a/K^o$ has connected structure group and is hence
orientable and as $Q_a/K^o=U'_a/G^o$ is a complex manifold and in particular
orientable we conclude that $Q_a$ is a compact orientable manifold.

We now apply duality to $U'_a$ and conclude that $H_i(U'_a,\Z)=0$ for
$2d-L<i<2d-h$, where $d$ is the (complex) dimension of $U'_a$ and $h:=\dim
H$. As $Q_a$ is homotopic to $U'_a$ the same is true for $Q_a$. Being an
orientable manifold we can apply duality to it together with the fact that $\dim
U'_a-\dim Q_a=\dim H-\dim K=h$. This gives $H^i(Q_a,\Z)=0$ for
$0<i<L-h$ and as $U'_a$ is connected so is $Q_a$ and we get that
$H^0(Q_a,\Z)=\Z$. Therefore we get that $H^i(U'_a,\Z)$ equals $\Z$ if $i=0$ and $0$
if $0<i<L-h$. We have a commutative diagram of $H$-torsors
\begin{displaymath}
\begin{CD}
U'_a @>>> \Eclass H\\
@VVV  @VVV\\
Y_a @>>> \Bclass H
\end{CD}
\end{displaymath}
and the map $U'_a \to \Eclass H$ induces an isomorphism in cohomological
degrees $<L-h$. This implies that the same is true for $Y_a \to
\Bclass H$ (by for instance the spectral sequence $H^i(H^j\times X)\implies
H^{i+j}(Y)$ that exists for any $H$-torsor $X \to Y$) and as $L$ can be made
arbitrarily large we are finished with the first part.

For the last part, we have by construction that $P/\mul=X_a$ and thus $P \to P/\mul$
is the $\mul$-torsor associated to the restriction $\sL$ of $\sO(1)$ to $X_a$ and
as $G$ acts compatibly on $P$ and $P/\mul$, $\sL$ extends to a line bundle $\sM$
on $P/H=Y_a$. Furthermore, the space of sections of $\sM^{\tensor r}$ equals
$(S^r V)^G$ which shows that $\sM^{\tensor k}$ is the restriction of the
$\sO(1)$ of $\P((S^k V)^G)$ to $Y_a \subseteq \P(V)//G \subseteq \P((S^k V)^G)$. 
\end{proof}
\end{theorem}
\begin{remark}
\part If $H$ is connected its centre is positive-dimensional precisely when $H$
is not semi-simple.

\part For a general reductive group $G$, $G\times \mul$ fulfils the conditions
of the theorem.

\part The proof of the theorem shows that we can refine (\ref{Stable topology}) to
$H^i_c(V_{ss}^n,\Z)=0$ for all $0< i\le K$ and $i \ne \dim G+1$ and $H^{\dim
G+1}_c(V_{ss}^n,\Z)=\Z$. It would be interesting to have a direct proof of this
fact.

\part Note that the appropriate homotopy version of the Lefschetz hyperplane
theorem (cf., \cite[Thm 7.4]{milnor63::morse}) gives that $X_a$ is a homotopy
approximation to to $\P(V)$. However, we need to have information on the
relation between $\P(V)\setminus T$ and $X_a\setminus T$. Our proof of such a
relation involved the use of duality. Duality is an intrinsically stable (in the
sense of homotopy theory) notion and hence destroys information on for instance
the fundamental group. Note incidentally that $T$ may very well have a larger
dimension than $X_a\setminus T$.

However, it seems reasonable to believe that the restriction to homological
approximation is just an artifice of the proof and I conjecture that there are
$G$-torsors over smooth and projective varieties such that the classifying map
induces an isomorphism on homotopy groups up to any limit.

\part One can, if one is so inclined, make an \emph{a posteriori} use of the
Lefschetz theorem to obtain a stronger result: Given an integer $L$ one can find
an algebraic $H$-torsor $P \to Y$ with $Y$ smooth projective of dimension $L$
such that $H^i(\Bclass H,\Z) \to H^i(Y,\Z)$ is an isomorphism for $i<L$ and
injective with a torsion free cokernel for $i=L$. Then $H^i(Y,\Z)$ for $i>L$ is
determined by duality. One constructs $Y$ from the $X$ of the theorem (using a
$K=L+1$) by taking the appropriate number of hyperplane sections.
\end{remark}
%;;%%Subsection Algebraic version
\begin{subsection}{Algebraic version}
\label{ss:Algebraic version}

In this subsection we shall quickly indicate the modifications needed to obtain
a theorem valid over an arbitrary (algebraically closed, for simplicity) field
$\k$, replacing of course classical cohomology with étale cohomology with
values in $\coeffs$, the inverse limit $\ili_n\Z/n\Z$, where $n$ runs over
integers invertible in $\k$.

The classifying map $X \to \Bclass H$ of an $H$-torsor \map{\pi}{P}{X} makes
sense if we replace the classifying space with the classifying algebraic
stack. Apart from that, the only non-algebraic component of the proof is the use
of the maximal compact subgroup $K$ and the reduction of the structure group of
a $H$-torsor to $K$. However, we only need the cohomological consequences of
that fact. More precisely we need to know that the dual (when the base is
smooth) of $R\pi_*\coeffs$ is isomorphic to $R\pi_*\coeffs[-h]$, where $h:=\dim
H$. For this we first note that the proof of the theorem shows that
$R^i\pi_*\coeffs$ is a constant sheaf of value $H^i(H,\coeffs)$. (Formally, the
argument shows this only for direct images with compact support but we can use
duality or smooth base change instead.) Suppose now that we know that
$H^i(H,\coeffs)=0$ for $i>h$, $H^h(H,\coeffs)=\coeffs$ and that the
multiplication map $R\Gamma(H,\coeffs)\Tensor^L R\Gamma(H,\coeffs) \to
R\Gamma(H,\coeffs)$ composed with $R\Gamma(H,\coeffs) \to H^h(H,\coeffs)[-h]$
gives an isomorphism $R\Gamma(H,\coeffs) \riso
R\Hom_{\coeffs}(R\Gamma(H,\coeffs),H^h(H,\coeffs))[-h]$. We then get a morphism
$R\pi_*\coeffs \to R^h\pi_*\coeffs[-h]=H^h(H,\coeffs)[-h]$ and it together with
the multiplication map $R\pi_*\coeffs\Tensor^L R\pi_*\coeffs \to R\pi_*\coeffs$
gives a morphism $R\pi_*\coeffs \to
R\sHom_{\coeffs}(R\pi_*\coeffs,H^h(H,\coeffs))[-h]$. To check that this is
isomorphism we can pull back by $\pi$ and by smooth base change the map becomes
just the constant map $R\Gamma(H,\coeffs) \riso
R\Hom_{\coeffs}(R\Gamma(H,\coeffs),H^h(H,\coeffs))[-h]$.

Now, in order to show duality for $R\Gamma(H,\coeffs)$ we start by noticing that
the cohomology of $H$ is isomorphic to that of $H/U$, where $U$ is the unipotent
radical of a Borel group $B$ of $H$. The quotient map $G/U \to G/B$ is a
$T$-torsor and by the argument just given and the fact $G/B$ is proper we are
reduced to proving duality for $R\Gamma(T,\coeffs)$. This of course is
well-known. (For the proof of duality for $R\Gamma(H,\coeffs)$ we could also
have started off from the truth of it over $\C$ and then used a specialisation
argument as in \cite[Sommes trig.:8.2]{deligne77::cohom} to show that the
cohomology is independent of the algebraically chosen base field and in
particular of the characteristic.)

In characteristic zero this gives us a purely algebraic proof of the theorem. In
positive characteristic the problem is that we do not know that a general member
of the linear system of $H$-invariant forms on $\P(V)$ is smooth outside of the
base (i.e., unstable) locus. This can be rectified by going directly to (in the
notations of the theorem) $X_a$: Letting $V_i:=\set{x \in \P(V)}{f_i(x)\ne 0}$
we have that the $V_i$ as well as their intersections are smooth and affine and
their union, $V'$, is the complement of $X_a$. Using the \v Cech spectral
sequence for $V'=\cup_iV_i$ we get that $H^i(V',\coeffs)=0$ for $i>\dim V'+a$
and using duality (as $V'$ is smooth) and the long exact sequence of cohomology
with compact support we get the needed properties of the cohomology of $X_a$.

Finally, as we no longer (in positive characteristic) can assume that
$X_a\setminus T$ is smooth we can not use it to get the smoothness for
$X_a/H$. However, $X_a/H$ is a complete intersection with respect to a very
ample linear system of $\P(V)//G$ lying in the smooth locus of $\P(V)//G$ so
that if the $f_i$ are chosen to be general the smoothness of $X_a/H$ follows
directly.

We can go even further and start with an algebraic group over an arbitrary (i.e.,
non-algebraically closed) field. We can then choose
our representation $V$ to be defined over the same field. The genericity
conditions for the choice of the $H$-invariant forms can be fulfilled over the
base field if it is infinite. In the case of a finite base field it can be
fulfilled by increasing the degree of the forms.
\end{subsection}
\end{section}
%;;%Section Applications
\begin{section}{Applications}

An example of how the theorem can be applied to finite groups is given in
\cite{Ek86:2}. There one uses the fact that for every prime $p$ there is a
finite group $G$ with elements $x,y,z \in H^1(G,\Z/p\Z)$ such that the Massey
product $\langle x,y,z\rangle$ is defined and non-trivial. Using a smooth and
projective approximation one concludes that there is a smooth and projective
variety with non-trivial $\Z/p\Z$-Massey products (on the contrary rational
Massey products of smooth and proper varieties are always trivial, cf.,
\cite{deligne75::real+kaehl}).

A somewhat vaguer application is to the theory of characteristic classes of
$G$-torsors. The theorem implies that a characteristic class of $G$-torsors is
determined by its restriction to the category of smooth projective varieties and
algebraic $G$-torsors for $G$ a reductive group.

Recall (cf., \cite{58::semin+c}) that an algebraic group $G$ (over an
algebraically closed field) is said to be \Definition{special} if every
algebraic $G$-torsor is trivial in the Zariski topology. Special groups have the
property that for every (topological) $G$-torsor $P \to X$, the associated
$G/B$-fibration $G/B\times_GP=:Y \to X$ has $H^*(Y)$ isomorphic to
$H^*(G/B\times X,\Z)$. We can use the theorem to prove a very specific version
of the converse of this result.
%;;%%Proposition Non-split
\begin{proposition}\label{Non-split}
Let $G$ be a connected non-special linear group. There is a smooth projective
variety $X$ and an algebraic $G$-torsor $P \to X$ such that if $Y \to X$ is the
associated $G/B$-fibration, where $B$ is a Borel subgroup of $G$, then $H^*(Y)$
is not isomorphic to $H^*(G/B\times X,\Z)$.
\begin{proof}
Suppose that we have proved the proposition instead for $G\times\mul$ (which
also is non-special) so that we have a $G\times\mul$-torsor $P \to X$ with the
requisite properties. If $B$ is a Borel subgroup of $G$, then $B\times\mul$ is a
Borel subgroup of $G\times\mul$. Hence the
$(G\times\mul)/(B\times\mul)$-fibration associated associated to $P \to X$ is
also the $G/B$-fibration associated to the $G$-torsor $P/\mul \to X$ and we get
the result for $G$.

Hence we may assume that $G$ has a centrally embedded $\mul$. We let $X$ be a
high degree approximation to $\Bclass G$ as in Theorem \ref{Approximation}. Then
the induced morphism $Y \to G/B\times_GEG$ is also highly (cohomologically)
connected. Now, $G/B\times_GEG$ is isomorphic to $EG/B=\Bclass(B)\sim \Bclass
T$, where $T$ is a maximal torus of $G$. This means that
$H^*(G/B\times_GEG,\Z)=H^*(\Bclass T,\Z)$ and $H^*(\Bclass T,\Z)$ is the
symmetric algebra on the cocharacter group of $T$. In particular it is torsion
free so if $H^*(\Bclass G,\Z)$ has torsion then we do not have an isomorphism
$H^*(\Bclass T)\iso H^*(G/B\times\Bclass T,\Z)$. As $X$ can be chosen to be an
arbitrarily high degree approximation if $H^*(\Bclass G,\Z)$ has torsion we can
find an $X$ fulfilling the conditions of the proposition. However, by
\cite[Exp~5: Thm 3,Thm 4]{58::semin+c} $H^*(\Bclass G,\Z)$ has torsion when $G$
is non-special.
\end{proof}
\end{proposition}
\begin{remark}
\part In many cases one can remove the condition that $G$ be reductive. For instance,
in characteristic zero any linear connected algebraic group $G$ contains a
complement to its unipotent radical $U$, i.e., a subgroup $H \subseteq G$ which
maps isomorphically to $G/U$. We can then construct an $H$-torsor is in the
proposition and extend it to a $G$-torsor through the inclusion $H \subseteq
G$. We then have that the associated $G/B$-fibration is isomorphic to the
$H/B'$-fibration associated to the original torsor.

\part In some cases at least one can prove the proposition through the
Atiyah-Hirzebruch result. For instance for $\SO_n$ one can look at the
diagonally embedded $(\Z/2)^{n-1}\subseteq\SO_n$. Its classifying space detects
all the torsion of $H^*(\Bclass\SO_n,\Z)$ so getting a
$(\Z/2)^{n-1}\times\mul$-torsor will give the proposition. Similarly, one can
use the $\Z/n$ Heisenberg group and its embedding into $\SL_n(\C)$ giving an
embedding of $(\Z/n)^2\subseteq\PSL_n(\C)$ which detects the torsion class of
$H^3(\PSL_n(\C),\Z)$. In characteristic $p$ and for a reductive group $G$
defined over a finite field one can use $G(\F)$ for a large enough finite
field. This inclusion does not necessarily lift to a homomorphism $G(\F) \to
G(\C)$ (in fact typically it does not) and only gives a map from $\Bclass G(\F)$
to the profinite completion (away from $p$) of $\Bclass G(\C)$ which is not
enough for our purposes (we would really need a group homomorphism $G(F) \to
G(\C)$).

\part Using the refinement discussed in the fifth part of the remark after
Theorem \ref{Approximation} one can get a more precise result in that one can
also demand that $H^*(Y,\Z)$ is torsion free and as $H^*(X,\Z)$ contains torsion
it is clear that $H^*(Y,\Z)$ is not isomorphic to $H^*(X\times G/B,\Z)$.
\end{remark}
Though it seems unlikely that the approximation theorem can be used to obtain
any substantial properties of $\Bclass G$ one can apply it to reprove a result
of Deligne (see \cite[Thm 9.1.1]{deligne74::theor+hodge+iii}) on the rational
cohomology of $\Bclass G$. Recall that $\Bclass G$ can be realised as (the
geometric realisation of) a simplicial scheme and hence its rational cohomology
has a natural mixed Hodge structure. The result of Deligne implies that the
mixed Hodge structure on $H^i(\Bclass G,\Q)$ is pure of weight $i$. This however
follows directly from the approximation theorem and the fact that the $i$'th
cohomology group of a smooth and proper variety is pure of weight $i$. Deligne's
result is more precise in that he proves that $H^i(\Bclass G,\Q)$ is zero for
odd $i$ and purely of type $(i/2,i/2)$ for even $i$. This, however, follows from
the weaker result: We have the Leray spectral sequence $H^i(\Bclass G,\Q)\Tensor
H^j(G/B,\Q) \implies H^{i+j}(\Bclass B,\Q)$ which is a spectral sequence of
mixed Hodge structures whose $E_2$-term is pure (of weight equal to the total
degree). Hence, the spectral sequence degenerates for weight reasons and we get
as a consequence the, \emph{splitting principle}, that $H^*(\Bclass G,\Q) \to
H^*(\Bclass B,\Q)$ is injective which reduces the stronger statement to a
computation of of the mixed Hodge structure on $H^*(\Bclass B,\Q)=H^*(\Bclass
T,\Q)$. Deligne goes in the other direction, using the degeneration of the
spectral sequence to get the splitting principle and hence the purity of
$H^*(\Bclass G,\Q)$. The degeneration of the spectral sequence is also
equivalent to the fact that $H^*(\Bclass T,\Q) \to H^*(G/B,\Q)$ is surjective
which was also proved by Grothendieck (see \cite[Exp~5: Cor 4, Lemme
10]{58::semin+c}). In any case the approximation theorem gives arguably a quite
natural explanation for the fact that $H^*(\Bclass G,\Q)$ is pure.

For our main application of the theorem we need to recall some facts. For a
field $\k$ one defines $\Kscheme{\k}$ as the group generated by elements $[X]$,
$X$ an algebraic $\k$-variety, such that isomorphic varieties give the same
element and $[X]-[Y]=[X\setminus Y]$ for a closed subvariety $Y \subseteq X$. It
has a ring structure characterised by $[X][Y]=[X\times Y]$ and one defines a
dimension filtration on the localisation $\Kscheme{\k}[\Lclass^{-1}]$, with
$\Lclass:=[\A^1]$, where $\text{Fil}^n\Kscheme{\k}[\Lclass^{-1}]$ is generated by
elements of the form $[X]\Lclass^{-i}$ for which $\dim X \le n+i$. Then
$\Kschemecompl{\k}$ is the completion of $\Kscheme{\k}[\Lclass^{-1}]$ with
respect to this filtration. If $G$ is a special algebraic group it follows
easily that for any algebraic $G$-torsor $P \to X$ we have $[P]=[G][X]$ in
$\Kscheme{\k}$. We aim to use Proposition \ref{Non-split} to give a (partial)
converse to this fact.
%;;%%Theorem
\begin{theorem}
For every connected non-special linear algebraic $G$ (over the field of complex
numbers) there is an algebraic $G$-torsor $P \to X$ such that $[P]\ne[G][X] \in
\Kschemecompl{\k}$.
\begin{proof}
As $G$ is non-special and its unipotent radical $U$ is special (cf.,
\cite[Exp~1, Prop~14]{58::semin+c}) we get that $G/U$ is non-special as
extensions of special groups are special (cf., \cite[Exp~1,
Lemme~6]{58::semin+c}) and we may assume that $G$ is reductive. In \cite[Prop.\
3.2]{ekedahl08::class} an invariant of elements of $\Kschemecompl{\k}$ is
defined whose value on $[X]$ for $X$ smooth and proper gives information
equivalent to the isomorphism class of $H^*(X,\Z)$ as a graded group. It is also
proved (cf., \cite[Lemma 3.8]{ekedahl08::class}) that if $[P]=[G][X]$ in
$\Kschemecompl{\k}$ for all $G$-torsors, then $[Y]=[G/B][X]=[G/B\times X]$ for
all $G/B$-fibrations $Y \to X$ with structure group $G$. However, if we had
$[Y]=[G/B\times X]$ for a smooth and proper $X$, then by what was just recalled
(and the fact that $G/B$ is smooth and proper) we would get $H^*(Y,\Z)\iso
H^*(G/B\times X,\Z)$ which by Proposition \ref{Non-split} is not always true.
\end{proof}
\end{theorem}
\begin{remark}
\part It is not difficult to extend the theorem to non-connected groups (which are
always non-special). In that case however $[P]=[G][B]$ should essentially never
be true, for instance if $G$ is finite it is true only if the torsor is trivial.

\part Note that having a base field of characteristic zero is (currently)
necessary: The invariant of \cite{ekedahl08::class} is defined using resolution
of singularities.

\part That invariant is additive and takes value in a torsion free abelian group
so that we get not only that $[P]\ne [G][X]$ but also that the difference
$[P]-[G][X]$ has infinite order.
\end{remark}
\end{section}
\bibliography{preamble,abbrevs,alggeom,algebra,topology,ekedahl}
\bibliographystyle{pretex}
\end{document}